\documentclass[12pt]{article}
\usepackage{fullpage}
\usepackage{amssymb,amsmath,amsfonts,amsthm,latexsym,enumerate,url,cases}
\numberwithin{equation}{section}
\usepackage{hyperref}
\hypersetup{colorlinks=true,citecolor=blue,linkcolor=blue,urlcolor=blue}

\newtheorem{theorem}{Theorem}[section] %
\newtheorem{lemma}[theorem]{Lemma} %
\usepackage[left=50pt,right=50pt,top=30pt,bottom=80pt]{geometry}
\begin{document}
\title{Additive representation functions and discrete convolutions}
\author{{Csaba S\'{a}ndor\footnote{csandor@math.bme.hu. The author was supported by the NKFIH Grant No. K129335 and by the Lend\" ulet program
of the Hungarian Academy of Sciences (MTA).}} \\
\small  $\dagger$Institute of Mathematics, Budapest University of Technology and Economics,\\
\small MTA-BME Lend\" ulet
Arithmetic Combinatorics Research Group\\
\small  H-1529 B. O. Box, Hungary
}
\date{}

\maketitle \baselineskip 18pt \maketitle \baselineskip 18pt

\begin{center}
\emph{This paper is dedicated to Professor Andr\' as S\' ark\" ozy on the occasion of his 80th birthday. }
\end{center}

{\bf Abstract.}
For a set $A$ of non-negative integers, let $R_A(n)$ denote the number of solutions to the equation $n=a+a'$ with $a$, $a'\in A$. Denote by $\chi_A(n)$ the characteristic function of $A$. Let $b_n>0$ be a sequence satisfying $\limsup_{n\to \infty}b_n<1$. In this paper, we prove some Erd\H os--Fuchs-type theorems about the error terms appearing in approximation formul\ae\ for $R_A(n)=\sum_{k=0}^n\chi _A(k)\chi _A(n-k)$ and $\sum_{n=0}^NR_A(n)$ having principal terms $\sum_{k=0}^nb_kb_{n-k}$ and $\sum_{n=0}^N\sum_{k=0}^nb_kb_{n-k}$, respectively.

 \vskip 3mm
 {\bf 2010 Mathematics Subject Classification:} 11B34

 {\bf Keywords and phrases:} Additive representation functions, Erd\H os--Fuchs theorem
\vskip 5mm

\section{Introduction}
Let $\mathbb{N}$ be the set of non-negative integers. For a subset $A\subseteq\mathbb{N}$, we define the representation function $R_A(n)$ of $A$ as the number of solutions to $n=a+a'$ with $a,a'\in A$. Let us denote the characteristic function of $A$ by $\chi_A(n)$; that is, $\chi_A(n)=1$ for $n\in A$, and $\chi_A(n)=0$ for $n\notin A$. Clearly, we have $R_A(n)=\sum_{k=0}^n\chi _A(k)\chi _A(n-k)$. Given now a sequence $b_n$ of real numbers with $0\le b_n\le 1$, we may be interested in comparing $R_A(n)$ and the accumulated representation function $\sum_{n=0}^NR_A(n)$ to the discrete convolution $\sum_{k=0}^nb_kb_{n-k}$ and the sum $\sum_{n=0}^N\!\sum_{k=0}^nb_kb_{n-k}$, respectively. A question that naturally arises is then: what can be said about the (asymptotic) error terms appearing? To address this question, let us begin by recalling the celebrated Erd\H os--Fuchs theorem, shown in \cite{EF}.{\parfillskip=60.9pt\par}

\begin{theorem}[Erd\H os--Fuchs, 1956]
  For any $c_1>0$, there is no subset $A\subseteq \mathbb{N}$ for which the relation $$\belowdisplayskip=0pt\sum_{n=0}^NR_A(n)=c_1N+o(N^{1/4}\log ^{-1/2}N)$$ is satisfied.
\end{theorem}

\vskip 2mm
In \cite{B}, Bateman extended this theorem as follows.
\begin{theorem}[Bateman, 1977]
Let $G(n)$ be a real-valued function on the non-negative integers that meets the following assumptions.
\begin{enumerate}
  \item  $G(n)\to \infty$ as $n\to \infty$.
  \item $\Delta^2G(n):=G(n)-2G(n-1)+G(n-2)\ge 0$ for $n$ sufficiently large.
  \item $\displaystyle\frac{G(2n)}{G(n)}<C\mkern2mu$ for some $C$ and for $n$ sufficiently large.
  \item $G(n)=o\big(\frac{n^2}{\log ^2n}\big)$
\end{enumerate}
Then, there does not exist any subset $A\subseteq\mathbb{N}$ for which the relation
$$\belowdisplayskip=0pt\mkern-10mu\sum_{n=0}^NR_A(n)=G(N)+o(G(N)^{1/4}\log ^{-1/2}N)$$ is satisfied.
\end{theorem}
\vskip 2mm

We may now formulate our first result, which is another extension of the Erd\H os--Fuchs Theorem.{\parfillskip=0pt\par}

\begin{theorem}\label{thm1}
Let $b_n$, $e_n$ be sequences of real numbers subject to the following conditions.
\begin{enumerate}
\item $b_n>0$ for $n$ sufficiently large.
\item $\,\displaystyle \frac{\sum_{k=0}^{2n}\,b_k}{\sum_{k=0}^n\,b_k}<C$ for some $C$ and for $n$ sufficiently large.
\item $\displaystyle \limsup_{n\to \infty} b_n\mkern-1mu<1$
\item $\displaystyle\mkern-3.5mu\lim_{n\to \infty} \frac{(\sum_{k=0}^nb_k^2)(\sum_{k=0}^ne_k^2)}{(\sum_{k=0}^nb_k)^3}=0$
\end{enumerate}
Then, there does not exist any subset $A\subseteq \mathbb{N}$ for which the relation
$$
\belowdisplayskip=0pt\sum_{n=0}^NR_A(n)=\sum_{n=0}^N\sum_{k=0}^nb_kb_{n-k} + e_N
$$
is satisfied.
\end{theorem}
Let us explore two particular cases of Theorem \ref{thm1}.

\pagebreak

First, with $c>0$, set $b_n=\sqrt{c}\frac{\binom{2n}{n}}{4^n}$ and $e_n\mkern-1mu=o\big(\frac{n^{1/4}}{\sqrt{\log n}}\big)$. Since $b_n\sim \frac{\sqrt{c}}{\sqrt{\pi}\sqrt{n} }$, we then have $\sum_{k=0}^nb_k^2=O(\log n)$ and $\sum_{k=0}^nb_k\mkern2mu\gg\sqrt{n}$. Also, $\mkern1mu\sum_{k=0}^ne_k^2\mkern1mu=o(\frac{n^{3/2}}{\log n})$. All the conditions in the theorem, including\parfillskip=0pt
$$
\belowdisplayskip=0pt
\displaystyle \lim_{n\to \infty} \frac{(\sum_{k=0}^nb_k^2)(\sum_{k=0}^ne_k^2)}{(\sum_{k=0}^nb_k)^3}=0,
$$\parfillskip=0pt plus 1fill
are hence met. Now, it is known that for $|z|<1$, we have $\frac{1}{(1-z)^{1/2}}=\sum_{n=0}^{\infty}\frac{\binom{2n}{n}}{4^n}z^n$; thus, we can write
$$
\sum_{n=0}^{\infty}cz^n=\left(  \frac{\sqrt{c}}{(1-z)^{1/2}} \right) ^2=\sum_{n=0}^{\infty }\bigg(\sum_{k=0}^nb_kb_{n-k}\bigg)z^n.
$$
It follows that
$$
\abovedisplayshortskip=-12pt
\mkern-6mu\sum_{n=0}^N\sum_{k=0}^nb_kb_{n-k}=c(N+1),
$$
and so Theorem \ref{thm1} implies the Erd\H os--Fuchs Theorem, indeed.

To see another interesting case, with $0<c_2<0.5$, set $b_0=\sqrt{2c_0}+\frac{c_1-3c_2}{\sqrt{2c_2}}$, $b_n=\sqrt{2c_2}$ for $n\ge 1$ and $e_n=o(\sqrt{n})$. Again, all the conditions in Theorem \ref{thm1}, including
$$
\displaystyle \lim_{n\to \infty} \frac{(\sum_{k=0}^nb_k^2)(\sum_{k=0}^ne_k^2)}{(\sum_{k=0}^nb_k)^3}=0,
$$
are met. Note that this time we have
$$
\sum_{n=0}^N\sum_{k=0}^nb_kb_{n-k}=c_2N^2+c_1N+O(1),
$$
allowing us to conclude the following.

\begin{theorem}\label{thm2}
Picking $0<c_2<0.5$, there does not exist any subset $A\subseteq \mathbb{N}$ for which the relation\parfillskip=0pt
$$
\belowdisplayskip=0pt\sum_{n=0}^NR_A(n)=c_2N^2+c_1N+o(\sqrt{N})
$$\parfillskip=0pt plus 1fill
is satisfied.
\end{theorem}

\vskip 2mm
In \cite{R}, Ruzsa shows that the error term appearing in the Erd\H os--Fuchs Theorem is almost sharp.{\parfillskip=0pt\par}%
\begin{theorem}[Ruzsa, 1997]
  There exists a subset $A\subseteq\mathbb{N}$ such that the relation
$$
\belowdisplayskip=0pt
\sum_{n=0}^NR_A(n)=\frac{\pi }{4}N+O(N^{1/4}\log N)
$$
is satisfied.

\end{theorem}

We will prove that the error term appearing in Theorem \ref{thm2} is almost sharp as well.

\begin{theorem}\label{thm3}
Picking integers $0<p<q$, there exists a subset $A\subseteq\mathbb{N}$ such that
$$
\belowdisplayskip=0pt
\sum_{n=0}^NR_A(n)=0.5\frac{p^2}{q^2}N^2+1.5\frac{p^2}{q^2}N+O(\sqrt{N}\log ^{1/2}N).
$$
\end{theorem}

Note that with $A=\mathbb{N}$, we have
$$
\sum_{n=0}^NR_{\mathbb{N}}(n)=0.5N^2+1.5N+1.
$$
This then shows that the assumption $c_2<0.5$ in Theorem \ref{thm2} cannot be relaxed.

Let us continue by recalling a result of Erd\H os and S\' ark\" ozy shown in \cite{ES1}.

\begin{theorem}[Erd\H os, S\' ark\" ozy, 1986]
  Let $F(n)$ be an arithmetic function such that $F(n)\to\infty$ as $n\to\mkern-1mu\infty$, $F(n)\le F(n+1)$ for $n\ge n_0$ and such that $F(n)=\mkern1mu o\big(\frac{n}{\log^2 n}\big)$. Then,
  $$
  \max _{0\leq n\le N}|F(n)-R_A(n)|=o\big(\sqrt{F(N)}\big)
  $$
  cannot hold for any subset $A\subseteq\mathbb{N}$.
\end{theorem}

The same authors show in \cite{ES2} that the error term above is almost sharp.

\begin{theorem}[Erd\H os, S\' ark\" ozy, 1986]\label{erdos-sarkozy}%
Let $F(n)$ be an arithmetic function satisfying $F(n)>36\log n$ for $n>n_0$. Picking an integer $n_1$, let $g(x)$ be a real-valued function defined on the (open) interval $(0,\infty)\mkern-2mu$ and $x_0\mkern-2mu$ be a real number subject to the following conditions.
\begin{enumerate}
  \item $g'(x)$ exists and is continuous on $(0,\infty)$.
  \item $g'(x)\le 0$ for $x\ge x_0$.
  \item $0<g(x)<1$ for $x\ge x_0$.
  \item $\left|F(n)-2{\displaystyle\int_0^{n/2}}g(x)g(n-x)dx\right|<(F(n)\log n)^{1/2}$ for $n>n_1$.
\end{enumerate}
Then, there exists a subset $A\subseteq\mathbb{N}$ such that $|R_A(n)-F(n)|<8(F(n)\log n)^{1/2}$ for $n\ge n_2$.
\end{theorem}

We will prove some analogous theorems.

\begin{theorem}\label{thm4}
Consider a sequence $b_n\mkern-3mu$ of real numbers subject to the following conditions.
\begin{enumerate}
\item $0\le b_n\le 1$
\item $\limsup_{n\to \infty}b_n<1$
\item $\sum_{k=0}^nb_kb_{n-k}\to \infty$ as $n\to \infty$.
\end{enumerate}
Then, there does not exist any subset $A\subseteq\mathbb{N}$ for which the relation
$$
\belowdisplayskip=0pt
R_A(n)=\sum_{k=0}^nb_kb_{n-k}+o\left(\bigg(\sum_{k=0}^nb_kb_{n-k}\bigg)^{1/2}\right)
$$
is satisfied.
\end{theorem}

Letting $b_n\mkern2mu=\sqrt{c}$ with $0<c<1$ yields the following result.
\begin{theorem}\label{thm5}
Picking $0<c<1$, there does not exist any subset $A\subseteq\mathbb{N}$ for which the relation
$$
\belowdisplayskip=0pt
R_A(n)=cn+o\big(\sqrt{n}\big)
$$
is satisfied.
\end{theorem}

Note that Theorem \ref{thm2} follows as a simple corollary to Theorem \ref{thm5}.

In the vein of Theorem \ref{erdos-sarkozy}, we will conclude by demonstrating that the error term appearing in Theorem \ref{thm4} is almost sharp as well. The proof will be based on a probabilistic argument.

\begin{theorem}\label{thm6}
Let $b_n$ be a sequence of real numbers with $0\le b_n\le 1$ such that $\sum_{k=0}^nb_kb_{n-k}>3\log n$ for $n\ge n_3$. Then, there exists a subset $A\subseteq\mathbb{N}$ such that we have
$$
\left|R_A(n)-\sum_{k=0}^nb_kb_{n-k}\right|<5\bigg(\log n\sum_{k=0}^nb_kb_{n-k}\bigg)^{1/2}
$$
for $n\ge n_4$.
\end{theorem}

As above, by letting $b_n\mkern2mu=\sqrt{c}$ with $0<c<1$, we may state the following consequence.
\begin{theorem}\label{thm7}
Picking $0<c<1$, there exists a subset $A\subseteq\mathbb{N}$ such that
$$
R_A(n)=cn+O(\sqrt{n}\log ^{1/2}n).
$$
\end{theorem}

\section{Proofs}\label{proofs}
{\bf Proof of Theorem \ref{thm1} .} We prove by contradiction. Let $b_n$, $e_n$ be sequences of real numbers as in the theorem. Suppose that there exists a subset $A\subseteq \mathbb{N}$ for which
$$
\sum_{n=0}^NR_A(n)=\sum_{n=0}^N\sum_{k=0}^nb_kb_{n-k} + e_N.
$$
Consider the generating function $A(z)$ of $A$; that is, write $A(z)=\sum _{a\in A}z^a$. For $|z|=r<1$, we have $A^2(z)=\sum_{n=0}^{\infty}R_A(n)z^n$ and $\frac{1}{1-z}A^2(z)=\sum_{N=0}^{\infty}(\sum_{n=0}^NR_A(n))z^N$. Defining now $f(z)=\sum_{n=0}^{\infty }b_nz^n$, we can write $f^2(z)=\sum_{n=0}^{\infty}(\sum_{k=0}^nb_kb_{n-k})z^n$ and $\frac{1}{1-z}f^2(z)=\sum_{N=0}^{\infty}(\sum_{n=0}^N(\sum_{k=0}^nb_kb_{n-k}))z^N$, yielding
$$
\frac{1}{1-z}A^2(z)=\sum_{N=0}^{\infty}\bigg(\sum_{n=0}^NR_A(n)\bigg)z^N=\sum_{N=0}^{\infty}\bigg(\sum_{n=0}^N\sum_{k=0}^nb_kb_{n-k} + e_N\bigg)z^N=\frac{1}{1-z}f^2(z)+\sum_{n=0}^{\infty}e_nz^n.
$$
Multiplying both sides by $1-z$, we get
\begin{equation}\label{eq7}
A^2(z)=f^2(z)+(1-z)\sum_{n=0}^{\infty}e_nz^n.
\end{equation}
For $0<r<1$, substitution then yields
\begin{align*}
A^2(r)&=f^2(r)+(1-r)\sum_{n=0}^{\infty}e_nr^n=(1-r)\bigg(\frac{1}{1-r}f^2(r)+\sum_{n=0}^{\infty}e_nr^n\bigg)\\
&=(1-r)\bigg(\sum_{N=0}^{\infty}\bigg(\bigg(\sum_{n=0}^N\sum_{k=0}^nb_kb_{n-k}\bigg) + e_N\bigg)r^N \bigg).
\end{align*}
Let us show that
$$
\mkern6mu e_N=\mkern2mu o\bigg(\sum_{n=0}^N\sum_{k=0}^nb_kb_{n-k}\bigg).
$$
By way of contradiction, suppose there is a constant $c>0$ such that for infinitely many $N$, we have
$$
|e_N|>c\bigg(\sum_{n=0}^N\sum_{k=0}^nb_kb_{n-k}\bigg).
$$
Note that the sequence $e_n$ necessarily has a non-zero term. Indeed, should $e_n\equiv 0$, we have $A^2(z)=\sum_{n=0}^{\infty}R_A(n)z^n=\sum_{n=0}^{\infty}(\sum_{k=0}^nb_kb_{n-k})z^n=f^2(z)$ and so $A(z)=\sum_{n=0}^{\infty}\chi _A(n)z^n=\sum_{n=0}^{\infty}b_nz^n=f(z)$. It would then follow that $b_n=\chi _A(n)$ for $n\in \mathbb{N}$, violating conditions 1 and 3. In particular, we see that condition 4 implies $\sum_{n=0}^Nb_n\to \infty$ as $N\to \infty$. Recall that by condition 1, there should exist an integer $n_5$ such that $b_n>0$ for $n\ge n_5$. Making use of condition 2 as well, this now allows us to write
$$
\sum_{n=0}^N\sum_{k=0}^nb_kb_{n-k}\ge \frac{1}{2}(b_{n_0}+b_{n_0+1}+\dots +b_{\lfloor \frac{N}{2} \rfloor })^2 -\bigg(\sum_{k=0}^{n_5-1}|b_k|\bigg)\bigg(\sum_{k=0}^n|b_k|\bigg) \ge \frac{1}{3C^2} (b_0+b_1+\dots +b_N)^2
$$
for $N$ large enough. It then follows that
$$
e_N^2>\frac{c^2}{9C^4}\bigg(\sum_{n=0}^Nb_n\bigg)^4\mkern-4mu,
$$
which, however, contradicts condition 4.

Having this shown, one can see that
\begin{equation}\label{eq4}
A(r)=(1+o(1))f(r)
\end{equation}
as $r\to 1^-$. Now, by condition 2,
\begin{align*}
\bigg(\sum_{n=0}^Nb_n\bigg)^3\mkern-6mu&\le C^6 \bigg(\sum_{n=0}^{\lfloor \frac{N}{3} \rfloor }b_n\bigg)^3\\
&\leq C^6\bigg(\mkern-6mu\sum_{\substack{(i,j,k)\\i+j+k\le N}}b_ib_jb_k + \bigg(\sum_{n<n_5}|b_n|\bigg)\bigg(\sum_{n_5\le n\le N}b_n\bigg)^2\mkern-6mu+\bigg( \sum_{n<n_5}|b_n|\bigg)^2\bigg(\sum_{n_0\le n\le N}b_n\bigg)  \bigg)\\
&\leq C^6\mkern-11mu\sum_{\substack{(i,j,k)\\i+j+k\le N}}b_ib_jb_k +o\bigg( \bigg(\sum_{n=0}^Nb_n\bigg)^3 \bigg)
\end{align*}
as $N\to \infty$, and so
$$
\frac{1}{2C^6}\bigg(\sum_{n=0}^Nb_n\bigg)^3\le\mkern-8mu\sum_{\substack{(i,j,k)\\i+j+k\le N}}b_ib_jb_k
$$
for $N$ large enough. Clearly, we have
$$
\mkern-8mu\bigg(\sum_{n=0}^Ne_n^2\bigg)\bigg(\sum_{n=0}^Nb_n^2\bigg)\ge\mkern-8mu\sum_{\substack{(i,j)\\i+j\le N}}e_i^2b_j^2\mkern2mu.
$$
It then follows from condition 4 that
$$
\lim_{N\to \infty}\frac{\sum_{\mkern-10mu\substack{(i,j)\\i+j\le N}}e_i^2b_j^2}{\sum_{\mkern-15mu\substack{(i,j,k)\\i+j+k\le N}}b_ib_jb_k}=0,
$$
and hence (see \cite{PSz} p. 21, Exercise 88)
\begin{equation}\label{eq11}
\lim_{r\to 1^-}\frac{(\sum_{n=0}^{\infty}b_n^2r^{2n})(\sum_{n=0}^{\infty}e_n^2r^{2n})}{(\sum_{n=0}^{\infty}b_nr^{2n})^3}=0.
\end{equation}
Write $S(z)=1+z+z^2+\dots +z^{M-1}$. Multiplying $S(z)^2$ by the expression for $A^2(z)$ in (\ref{eq7}), we get
$$
A^2(z)S^2(z)=f(z)^2S(z)^2+(1-z^M)S(z)\sum_{n=0}^{\infty}e_nz^n,
$$
yielding
$$
|A(z)S(z)|^2\le M^2|f(z)|^2+2|S(z)|\bigg|\sum_{n=0}^{\infty}e_nz^n\bigg|.
$$
Integrating both sides of the above inequality with respect to $\varphi$ in $z=re^{i\varphi}$ leads to the inequality\parfillskip=0pt
\begin{equation}\label{eq9}
\int_0^{2\pi }|A(re^{i\varphi})S(re^{i\varphi})|^2d\varphi \le M^2\int_0^{2\pi}|f(re^{i\varphi})|^2d\varphi +2\int_0^{2\pi }|S(re^{i\varphi})|\bigg|\sum_{n=0}^{\infty}e_nr^ne^{in\varphi}\bigg|d\varphi.
\end{equation}\parfillskip=0pt plus 1fill
Write now $A(z)S(z)=\sum_{n=0}^{\infty}c_nz^n$ for some $c_n\in \mathbb{N}$. With $c_n$ being an integer, $c_n\le c_n^2$. By Parseval's theorem and equation (\ref{eq4}), we then get
\begin{align*}
\int_0^{2\pi} |A(re^{i\varphi})S(re^{i\varphi})|^2d\varphi&=2\pi \sum_{n=0}^{\infty} c_n^2r^{2n}\\
&\mkern-48mu\ge 2\pi \sum_{n=0}^{\infty} c_nr^{2n}=A(r^2)S(r^2)=2\pi (1+o(1))(1+r^2+r^4+\dots +r^{2M-2})\sum_{n=0}^{\infty}b_nr^{2n}
\end{align*}
for $0<r<1$. A repeated application of Parseval's theorem yields
$$
\int_0^{2\pi} |f(re^{i\varphi})|^2d\varphi =2\pi \sum_{n=0}^{\infty}b_n^2r^{2n}.
$$
By the Cauchy--Schwarz inequality and using Parseval's theorem once again, we can write
\begin{align*}
\mkern-42mu\int_0^{2\pi}|S(re^{i\varphi})|\bigg|\sum_{n=0}^{\infty}e_nr^ne^{in\varphi }\bigg|d\varphi&\le\Bigg(\int_0^{2\pi }|S(re^{i\varphi})|^2d\varphi \int_0^{2\pi}\bigg|\sum_{n=0}^{\infty}e_nr^ne^{in\varphi }\bigg|^2d\varphi\Bigg)^{1/2}\\[4pt]
&=\Bigg(2\pi \sum_{n=0}^{M-1}r^{2n}\cdot 2\pi \sum_{n=0}^{\infty}e_n^2r^{2n}\Bigg)^{1/2}\le 2\pi\Bigg(M\sum_{n=0}^{\infty}e_n^2r^{2n}\Bigg)^{1/2}.
\end{align*}
It then follows from  (\ref{eq11}) and (\ref{eq9}) that
{\abovedisplayskip=0pt\belowdisplayskip=0pt\begin{align*}
2\pi (1+o(1))(1+r^2+r^4+\dots +r^{2M-2})\sum_{n=0}^{\infty}b_nr^{2n}&\le2\pi M^2\sum_{n=0}^{\infty}b_n^2r^{2n} + 2\pi \Bigg(M\sum_{n=0}^{\infty}e_n^2r^{2n}\Bigg)^{1/2}\\
&\leq 2\pi M^2\sum_{n=0}^{\infty}b_n^2r^{2n} +2\pi  \Bigg(o(1)M\frac{(\sum_{n=0}^{\infty} b_nr^{2n})^3}{\sum_{n=0}^{\infty}b_n^2r^{2n}}\Bigg)^{1/2},
\end{align*}}
and hence
\begin{equation}\label{eq3}
(1+o(1))(1+r^2+r^4+\dots +r^{2M-2})\sum_{n=0}^{\infty}b_nr^{2n} \le M^2\sum_{n=0}^{\infty}b_n^2r^{2n} + o\!\left(\!\Bigg(M\frac{(\sum_{n=0}^{\infty} b_nr^{2n})^3}{\sum_{n=0}^{\infty}b_n^2r^{2n}}\Bigg)^{1/2}\right)
\end{equation}
as $r\to 1^-$. This leaves us with two cases to inspect.

$\mkern.8mu$Case 1. If $\limsup_{r\to 1^-}\frac{\sum_{n=0}^{\infty}b_nr^{2n}}{\sum_{n=0}^{\infty}b_n^2r^{2n}}=\infty$, then there must exist a strictly increasing sequence $r_k<1$ such that $\lim _{k\to \infty}r_k=1$ and $\lim_{k\to \infty} \frac{\sum_{n=0}^{\infty}b_nr_k^{2n}}{\sum_{n=0}^{\infty}b_n^2r_k^{2n}}=\infty$.
We claim that
\begin{equation}\label{eq6}
\mkern-4mu\lim_{k\to \infty} \frac{\sum_{n=0}^{\infty}b_nr_k^{2n}}{\sum_{n=0}^{\infty}b_n^2r_k^{2n}}=o\left(\frac{1}{1-r_k}\right).
\end{equation}
Indeed, if $\sum_{n=0}^{\infty}b_n^2=\infty$, then $\lim_{k\to \infty}\sum_{n=0}^{\infty}b_n^2r_k^{2n}=\infty$ and $\sum_{n=0}^{\infty}b_nr_k^{2n}=O(\sum_{n=0}^{\infty}r_k^{2n})=O\big(\frac{1}{1-r_k}\big)$, from which (\ref{eq6}) follows. If $\mkern1mu\sum_{n=0}^{\infty}b_n^2\mkern1mu<\infty$, then by the Cauchy--Schwarz inequality, for every positive integer $N$, we can write
$$
\abovedisplayskip=4pt
\bigg|\sum_{n=0}^Nb_nr_k^{2n}\bigg|\le \Bigg(\bigg(\sum_{n=0}^Nb_n^2\bigg)\bigg(\sum_{n=0}^Nr_k^{2n}\bigg)\Bigg)^{1/2}=O\!\left(\!\Bigg(\frac{1}{1-r_k}\Bigg)^{1/2}\right)=o\left(\frac{1}{1-r_k}\right),
$$
and the same conclusion can be drawn. Define now
$$
\advance\abovedisplayskip by 4pt
M=M_k=0.5\bigg\lfloor\frac{\sum_{n=0}^{\infty}b_nr_k^{2n}}{\sum_{n=0}^{\infty}b_n^2r_k^{2n}}\bigg\rfloor =o\left(\frac{1}{1-r_k}\right)
$$
so that we have
$$
\abovedisplayshortskip=-4pt
1+r_k^2+r_k^4+\dots +r_k^{2M-2}=(1+o(1))M
$$
as $k\to \infty$.
By (\ref{eq3}), we get
$$
\mkern-4mu\left(\frac{1}{2}\mkern-1mu+\mkern-1mu o(1)\!\right)\frac{\sum_{n=0}^{\infty}b_nr_k^{2n}}{\sum_{n=0}^{\infty}b_n^2r_k^{2n}}\sum_{n=0}^{\infty}b_nr_k^{2n} \le \frac{1}{4}\left( \frac{\sum_{n=0}^{\infty}b_nr_k^{2n}}{\sum_{n=0}^{\infty}b_n^2r_k^{2n}}\right)^{\!2}\mkern-4mu\sum_{n=0}^{\infty}b_n^2r_k^{2n}  + o\!\left(\mkern-4mu\Bigg(\frac{\sum_{n=0}^{\infty}b_nr_k^{2n}}{\sum_{n=0}^{\infty}b_n^2r_k^{2n}}\frac{(\sum_{n=0}^{\infty} b_nr_k^{2n})^3}{\sum_{n=0}^{\infty}b_n^2r_k^{2n}}\Bigg)^{\mkern-6mu 1/2}\mkern-1mu\right)
$$
as $k\to \infty$, a contradiction.

Case 2. If $\limsup_{r\to 1^-}\frac{\sum_{n=0}^{\infty}b_nr^n}{\sum_{n=0}^{\infty}b_n^2r^n}<\infty$, then setting $M=1$ in equation (\ref{eq3}) yields
$$
(1+o(1))\sum_{n=0}^{\infty}b_nr^{2n} \le \sum_{n=0}^{\infty}b_n^2r^{2n} + o\!\left(\!\Bigg(\frac{(\sum_{n=0}^{\infty} b_nr^{2n})^3}{\sum_{n=0}^{\infty}b_n^2r^{2n}}\Bigg)^{1/2}\right).
$$
Note that
$$
o\!\left(\!\Bigg(\frac{(\sum_{n=0}^{\infty} b_nr^{2n})^3}{\sum_{n=0}^{\infty}b_n^2r^{2n}}\Bigg)^{1/2}\right)=o\bigg(\sum_{n=0}^{\infty} b_nr^{2n}\bigg),
$$
so we get
$$
(1+o(1))\sum_{n=0}^{\infty}b_nr^{2n} \le \sum_{n=0}^{\infty}b_n^2r^{2n}
$$
as $r\to 1^-$, a contradiction again as $\lim_{r\to 1-}\sum_{n=0}^{\infty }b_nr^{2n}=\infty $ and $\limsup_{n\to \infty} b_n<1$.

\bigskip

{\bf Proof of Theorem \ref{thm3}.} For each $n\in \mathbb{N}$, consider a simple random sample $X_n$ of size $|X_n|=p$ with members drawn (uniformly and independently) from the set $\{ nq,\mkern2mu nq+1,\dots ,nq+q-1\}$. Let $A=\bigcup _{n=0}^{\infty} X_n\subseteq\mathbb{N}$. Our goal is to prove that $\sum_{n=0}^NR_A(n)=0.5\frac{p^2}{q^2}N^2+1.5\frac{p^2}{q^2}N+O(\sqrt{N}\log ^{1/2}N)$ holds with probability 1.

To this end, define $Y_{u,v}$ for $(u,v)\in\mathbb{N}\times\mathbb{N}$ as the number of pairs $(a,a')$ with $a,a'\in A$ that satisfy
$$
uq\le a\le uq+q-1,\text{ }vq\le a'\le vq+q-1\text{ and }a+a'\le N.
$$
On the one hand, note that if the inequalities
$$
\mkern16mu u+v\le \bigg\lfloor \frac{N}{q}\bigg\rfloor -2,\text{ }uq\le a\le uq+q-1\,\text{ and }\,vq\le a'\le vq+q-1
$$
hold, then $a+a'\le N$. On the other hand, if
$$
\mkern16mu u+v> \bigg\lfloor \frac{N}{q}\bigg\rfloor,\,\;\phantom{-2}\text{ }uq\le a\le uq+q-1\,\text{ and }\,vq\le a'\le vq+q-1
$$
hold, then $a+a'>N$.
Now, as per the definition of $R_A(n)$, we have
\begin{equation}\label{eq1.1}
\begin{split}
 \sum_{n=0}^NR_A(n)&=\sum_{M=0}^{\lfloor \frac{N}{q}\rfloor -2}\sum_{m=0}^M Y_{m,M-m}+
 \sum_{m=0}^{\lfloor \frac{N}{q}\rfloor -1}Y_{m,\lfloor \frac{N}{q}\rfloor -1-m}+
 \sum_{m=0}^{\lfloor \frac{N}{q}\rfloor }Y_{m,\lfloor \frac{N}{q}\rfloor -m}\\
&=\sum_{M=0}^{\lfloor \frac{N}{q}\rfloor -2}\sum_{m=0}^M p^2+
 \sum_{m=0}^{\lfloor \frac{N}{q}\rfloor -1}Y_{m,\lfloor \frac{N}{q}\rfloor -1-m}+
 \sum_{m=0}^{\lfloor \frac{N}{q}\rfloor }Y_{m,\lfloor \frac{N}{q}\rfloor -m},
\end{split}
\end{equation}
and hence the expected value $\mathbf{E}(\sum_{n=0}^NR_A(n))$ may be calculated as
\begin{equation}\label{eq1.2}
\begin{split}
 &\mathbf{E}\bigg(\sum_{n=0}^NR_A(n)\bigg)=\sum_{M=0}^{\lfloor \frac{N}{q}\rfloor -2}\sum_{m=0}^M p^2+
 \mathbf{E}\Bigg(\sum_{m=0}^{\lfloor \frac{N}{q}\rfloor -1}Y_{m,\lfloor \frac{N}{q}\rfloor -1-m}\Bigg)+
 \mathbf{E}\Bigg(\sum_{m=0}^{\lfloor \frac{N}{q}\rfloor }Y_{m,\lfloor \frac{N}{q}\rfloor -m}\Bigg).
\end{split}
\end{equation}
Let us introduce the sets $R_{u,v}=\{(x,y)\in\mathbb{N}\times\mathbb{N}\;|\;uq\le x\le uq+q-1,\,vq\le y\le vq+q-1\}$ and $S_N=\{(x,y)\in\mathbb{N}\times\mathbb{N}\;|\;x+y\le N\}$. It is easy to see that $\mathbb{P}(k\in A)=\frac{p}{q}$ for all $k\in \mathbb{N}$, so for all pairs $(k,l)\in\mathbb{N}\times\mathbb{N}$, we may conclude the following.
$$\mathbf{E}(k\in A,l\in A)=\begin{cases}
\mkern12mu\frac{p^2}{q^2}&\mbox{if }\,iq\le k\le iq+q-1\mbox{ and }jq\le k\le jq+q-1\mbox{ with }i\ne j\\
\frac{\binom{q-2}{p-2}}{\binom{q}{p}}=\frac{p(p-1)}{q(q-1)}&\mbox{if }\,iq\le k,l\le iq+q-1\mbox{ with }k\ne l
\end{cases}\mkern-4mu$$
Except for when $u=v=\big\lfloor \frac{N}{2q}\big\rfloor$, this then tells us that $\mathbf{E}(Y_{u,v})=\frac{p^2}{q^2}|R_{u,v}\cap S_N|$, and hence
\bgroup
\belowdisplayskip=0pt
\begin{align*}
\mkern-14mu\mathbf{E}\bigg(\sum_{n=0}^NR_A(n)\bigg)=\sum_{(u,v)}\mathbf{E}(Y_{u,v})&=O(1)+\sum_{(u,v)}\frac{p^2}{q^2}|R_{u,v}\cap S_N|\\
&=O(1)+\frac{p^2}{q^2}|S_N|\\
&=O(1)+\frac{p^2}{q^2}\frac{(N+1)(N+2)}{2}.
\end{align*}
\egroup
It now follows from (\ref{eq1.1}) and (\ref{eq1.2}) that
\begin{align*}
&\bigg|\sum_{n=0}^NR_A(n)-\bigg(0.5\frac{p^2}{q^2}N^2+1.5\frac{p^2}{q^2}N\bigg)\bigg|=O(1)+\bigg|\sum_{n=0}^NR_A(n)-\mathbf{E}\bigg(\sum_{n=0}^NR_A(n)\bigg)\bigg|\\
&=O(1)+\bigg|\sum_{m=0}^{\lfloor \frac{N}{q}\rfloor -1}Y_{m,\lfloor \frac{N}{q}\rfloor -1-m}\mkern2mu+\mkern2mu\sum_{m=0}^{\lfloor \frac{N}{q}\rfloor }\mkern4mu Y_{m,\lfloor \frac{N}{q}\rfloor -m}\mkern2mu-\bigg(\mathbf{E}\bigg(\sum_{m=0}^{\lfloor \frac{N}{q}\rfloor -1}Y_{m,\lfloor \frac{N}{q}\rfloor -1-m}\bigg)+\mathbf{E}\bigg(\sum_{m=0}^{\lfloor \frac{N}{q}\rfloor }Y_{m,\lfloor \frac{N}{q}\rfloor -m}\bigg)\mkern-4mu\bigg)\bigg|\\
&\le O(1)+\bigg|\sum_{m=0}^{\lfloor \frac{N}{q}\rfloor -1}Y_{m,\lfloor \frac{N}{q}\rfloor -1-m}-\mathbf{E}\bigg(\sum_{m=0}^{\lfloor \frac{N}{q}\rfloor -1}\mkern-4mu Y_{m,\lfloor \frac{N}{q}\rfloor -1-m}\bigg)\bigg|+\bigg|\sum_{m=0}^{\lfloor \frac{N}{q}\rfloor }\bigg(\mkern-2mu Y_{m,\lfloor \frac{N}{q}\rfloor -m}-\mathbf{E}\bigg(\sum_{m=0}^{\lfloor \frac{N}{q}\rfloor }Y_{m,\lfloor \frac{N}{q}\rfloor -m}\bigg)\mkern-4mu\bigg)\bigg|.
\end{align*}
At this point, we are to make use of Hoeffding's inequality (see for example \cite{H}).
\begin{lemma}[Hoeffding's inequality]
  Let $\xi _1,\xi _2,\ldots,\xi_k$ be independent random variables bounded by the intervals $[a_1,b_1],\ldots,[a_k,b_k]$, respectively; that is, $a_i\le \xi _i\le b_i$. Let $\sum_{i=1}^k(b_i-a_i)^2\le D^2$, and write $\eta =\sum_{i=1}^k\xi _i$. Then, for all $y\geq 0$, the inequality $\mathbf{P}(|\eta -\mathbf{E}(\eta)|\ge yD)\le exp(-2y^2)$ holds.
\end{lemma}
Applying Hoeffding's inequality to the random variables $Y_{0,\lfloor \frac{N}{q}\rfloor -1},Y_{1,\lfloor \frac{N}{q}\rfloor -2},\ldots,Y_{\lfloor \frac{N}{q}\rfloor -1,0}$ and to $Y_{0,\lfloor \frac{N}{q}\rfloor },Y_{1,\lfloor \frac{N}{q}\rfloor -1},\ldots ,Y_{\lfloor \frac{N}{q}\rfloor ,0}$ with $0\le Y_{u,v}\le p^2$ and $D^2=p^2N$, we get the upper bounds
$$
\belowdisplayskip=0pt
\mkern-44mu p_N\mkern2mu=\mkern2mu\mathbf{P}\Bigg(\bigg|\sum_{m=0}^{\lfloor \frac{N}{q}\rfloor -1}Y_{m,\lfloor \frac{N}{q}\rfloor -1-m}-\mathbf{E}\bigg(\sum_{m=0}^{\lfloor \frac{N}{q}\rfloor -1}Y_{m,\lfloor \frac{N}{q}\rfloor -1-m}\bigg)\bigg|>\sqrt{p^2N} \sqrt{\log N}\Bigg)\le \frac{1}{N^2}
$$
and
$$
\mkern-106mu p'_N\mkern2mu=\mkern2mu\mathbf{P}\Bigg(\bigg|\sum_{m=0}^{\lfloor \frac{N}{q}\rfloor }Y_{m,\lfloor \frac{N}{q}\rfloor -m}-\mathbf{E}\bigg(\sum_{m=0}^{\lfloor \frac{N}{q}\rfloor }Y_{m,\lfloor \frac{N}{q}\rfloor -m}\bigg)\bigg|>\sqrt{p^2N} \sqrt{\log N}\Bigg)\le \frac{1}{N^2}.
$$
Both $p_N$, $p'_N$ having a finite sum over $N\in\mathbb{N}$, the Borel--Cantelli lemma tells us that with probability 1, only finitely many of the corresponding events occur. Putting all this together, finally, we see that
$$
\sum_{n=0}^NR_A(n)=0.5\frac{p^2}{q^2}N^2+1.5\frac{p^2}{q^2}N+O(\sqrt{N}\log ^{1/2}N)
$$
holds with probability $1$, and so the proof is complete.

\bigskip

{\bf Proof of Theorem \ref{thm4}.}
By contradiction, let $b_n$ be a sequence of real numbers subject to the conditions stated, and suppose that we have $R_A(n)=\sum_{k=0}^nb_kb_{n-k}+e_n$ with $e_n=o(\sqrt{\sum_{k=0}^nb_kb_{n-k}})$. With notation as in the proof of Theorem \ref{thm1}, we can write
$$
A^2(z)=\sum_{n=0}^{\infty }R_A(n)z^n=\sum_{n=0}^{\infty}\bigg(\sum_{k=0}^nb_kb_{n-k}+e_n\bigg)z^n=f^2(r)+\sum_{n=0}^{\infty}e_nz^n\mkern-2mu,
$$
and since $\sum_{k=0}^nb_kb_{n-k} \to \infty$ as $n\to \infty$ and $e_n=o(\sum_{k=0}^nb_kb_{n-k})$, one can also see that
\begin{equation}\label{eq10}
A(r)=(1+o(1))f(r)
\end{equation}
as $r\to 1^-$. Integrating now both sides of the (triangle) inequality
$$
|A(z)|^2\le |f(z)|^2+\bigg|\sum_{n=0}^{\infty}e_nz^n\bigg|
$$
with respect to $\varphi$ in $z=re^{i\varphi}$, for $0<r<1$, we get the inequality
\begin{equation}\label{eq5}
\int _{0}^{2\pi } |A(re^{i\varphi })|^2d\varphi\le \int _{0}^{2\pi } |f(re^{i\varphi })|^2d\varphi+\int _{0}^{2\pi }\bigg|\sum_{n=0}^{\infty}e_nr^ne^{in\varphi}\bigg|d\varphi.
\end{equation}
Parseval's theorem and equation \ref{eq10} then show that
$$
\int _{0}^{2\pi }|A(re^{i\varphi})|^2d\varphi =2\pi A(r^2)=2\pi (1+o(1))f(r^2)=(1+o(1))2\pi \sum_{n=0}^{\infty}b_nr^{2n}
$$
and
$$
\abovedisplayshortskip=-12pt
\int_0^{2\pi}|f(re^{i\varphi })|^2d\varphi =2\pi \sum_{n=0}^{\infty}b_n^2r^{2n}$$
as $r\to 1^-$, so by the Cauchy--Schwarz inequality and applying Parseval's theorem once again, we get
\begin{align*}
\int_0^{2\pi }\bigg|\sum_{n=0}^{\infty}e_nr^ne^{2i\varphi n}\bigg|d\varphi&\le\left(\!\bigg(\int_0^{2\pi}1d\varphi\bigg)\bigg(\int_0^{2\pi}\bigg|\sum_{n=0}^{\infty}e_nr^ne^{i\varphi n}\bigg|^2d\varphi\bigg)\!\right)^{1/2}\mkern-2mu=\bigg(2\pi\cdot 2\pi \sum_{n=0}^{\infty}e_n^2r^{2n}\bigg)^{1/2}\\[6pt]
&=\left(2\pi\cdot o\bigg(\sum_{n=0}^{\infty}\bigg(\sum_{k\le n}b_kb_{n-k}\bigg)r^{2n}\bigg)\!\right)^{1/2}\!\!=o(\sqrt{f^2(r^2)})=o\bigg(\sum_{n=0}^{\infty}b_nr^{2n}\bigg)
\end{align*}
as $r\to 1^-$. Note, however, that the above inequality and equation (\ref{eq5}) now yield
$$
(1+o(1))2\pi \sum_{n=0}^{\infty}b_nr^{2n}\le 2\pi \sum_{n=0}^{\infty}b_n^2r^{2n}
$$
as $r\to 1^-$, a contradiction as $\limsup_{r\to 1-}\sum_{n=0}^{\infty }b_nr^{2n}=\infty $ and $\limsup _{n\to \infty}b_n<1$.

\bigskip

{\bf Proof of Theorem \ref{thm6}.} Consider a random subset $A\subseteq\mathbb{N}$ with $\mathbf{P}(n\in A)=b_n$ for $n\in \mathbb{N}$. We are to show that
$$
\abovedisplayshortskip=-6pt
\bigg|R_A(n)-\sum_{k=0}^nb_kb_{n-k}\bigg|<8\bigg(\!\log n\sum_{k=0}^nb_kb_{n-k}\bigg)^{1/2}
$$
holds for $n$ large enough with probability 1.

To this end, note that we can write
$$
\mkern-7mu\mathbf{E}(R_A(n))=2\mkern-4mu\sum_{i=0}^{\lceil \frac{n}{2}\rceil -1}\mkern-4mu\mathbf{E}(i\in A,\,n-i\in A)+\mkern-23mu\sum_{\lceil \frac{n}{2}\rceil -1<i<n-(\lceil \frac{n}{2}\rceil -1)}\mkern-23mu\mathbf{E}(i\in A,\,n-i\in A)=\sum_{k=0}^nb_kb_{n-k}+e_n
$$
with $|e_n|\le 1$. We now wish to deploy the following Chernoff bound variant (see Corollary 1.9 in \cite{TV}).
\begin{lemma}[Chernoff bound]\label{chernoff}%
  Let $t_1,\ldots,t_k$ be independent random variables taking boolean values. Consider $X=t_1+t_2+\dots +t_k$. Then, $\mathbf{P}(|X-\mathbf{E}(X)|\ge \varepsilon\, \mathbf{E}(X))\le 2e^{-\min (\varepsilon ^2/4,\;\varepsilon /2)\,\mathbf{E}(X)}$ for all $\varepsilon >0$.\parfillskip=0pt
\end{lemma}
For $0\le i\le \lceil \frac{n}{2}\rceil -1$, define $t_i$ as the indicator random variable of the event: $i\in A$ and $n-i\in A$. Clearly, we have
$$
\advance\abovedisplayshortskip by -10pt
\advance\belowdisplayshortskip by 2pt
\mathbf{E}\bigg(\sum_{i=0}^{\lceil \frac{n}{2} \rceil -1}t_i\bigg)=\sum_{i=0}^{\lceil \frac{n}{2}\rceil -1}b_ib_{n-i},
$$
and it is also easy to check that the inequality
$$
\advance\abovedisplayshortskip by-4pt
|R_A(n)-\mathbf{E}(R_A(n))|\ge 5\bigg(\log n\sum_{k=0}^nb_kb_{n-k}\bigg)^{1/2}
$$
implies
$$
\Bigg|\sum_{i=0}^{\lceil \frac{n}{2}\rceil -1} t_i\mkern4mu-\mkern-4mu\sum_{i=0}^{\lceil \frac{n}{2}\rceil -1}b_ib_{n-i}\Bigg|\ge\Bigg(\frac{8\log n}{\sum_{k=0}^{\lceil \frac{n}{2} \rceil -1}b_kb_{n-k}}\Bigg)^{1/2}\cdot\sum_{k=0}^{\lceil \frac{n}{2} \rceil -1}b_kb_{n-k}
$$
for $n$ large enough. By Lemma \ref{chernoff}, setting $\displaystyle\varepsilon =\Big(\frac{8\log n}{\sum_{k=0}^{\lceil \frac{n}{2}\rceil -1}b_kb_{n-k}}\Big)^{1/2}\!<2$, we get
\bgroup
\advance\belowdisplayskip by 4pt
\begin{align*}
&\mkern-66mu p_n\mkern2mu=\mkern2mu\mathbf{P}\left(|R_A(n)-\mathbf{E}(R_A(n))|\ge 5\bigg(\log n\sum_{k=0}^nb_kb_{n-k}\bigg)^{1/2} \right)\\
&\mkern-46mu\leq\mathbf{P}\left(\Bigg|\sum_{i=0}^{\lceil \frac{n}{2}\rceil -1} t_i\mkern4mu-\mkern-4mu\sum_{i=0}^{\lceil \frac{n}{2}\rceil -1}b_ib_{n-i}\Bigg|\ge\Bigg(\frac{8\log n}{\sum_{k=0}^{\lceil \frac{n}{2} \rceil -1}b_kb_{n-k}}\Bigg)^{1/2}\cdot\sum_{k=0}^{\lceil \frac{n}{2} \rceil -1} b_kb_{n-k} \right)\\[-2pt]
&\mkern-46mu\leq\exp\Bigg(-\frac{1}{4}\frac{8\log n}{\sum_{k=0}^{\lceil \frac{n}{2} \rceil -1}b_kb_{n-k}}\sum_{k=0}^{\lceil \frac{n}{2} \rceil -1}b_kb_{n-k}\Bigg)=e^{-2\log n}=\frac{1}{n^2}
\end{align*}
\egroup
for $n$ large enough. Since the probabilities $p_n$ have a finite sum over $n\in\mathbb{N}$, the Borel--Cantelli lemma tells us that with probability $1$, only finitely many of the corresponding events occur. All in all,
$$\bigg|R_A(n)-\sum_{k=0}^nb_kb_{n-k}\bigg|<8\bigg(\!\log n\sum_{k=0}^nb_kb_{n-k}\bigg)^{1/2}$$
holds for $n$ large enough with probability $1$, which then concludes our argument.

\bigskip

\renewcommand{\refname}{Bibliography}

\end{document}